# INERTIAL MANIFOLDS FOR REACTION-DIFFUSION EQUATIONS ON GENUINELY HIGH-DIMENSIONAL THIN DOMAINS

M. PRIZZI AND K. P. RYBAKOWSKI

ABSTRACT. In this paper we study a family of semilinear reaction-diffusion equations on spatial domains $\Omega_\epsilon$, $\epsilon > 0$, in $\mathbb{R}^\ell$ lying close to a $k$-dimensional submanifold $\mathcal{M}$ of $\mathbb{R}^\ell$. As $\epsilon \to 0^+$, the domains collapse onto (a subset of) $\mathcal{M}$. As it was proved in a previous paper (M. Prizzi, M. Rinaldi and K. P. Rybakowski, Curved thin domains and parabolic equations, Stud. Math. 151), the above family has a limit equation, which is an abstract semilinear parabolic equation defined on a certain limit phase space denoted by $H_s^1(\Omega)$. The definition of $H_s^1(\Omega)$, given in the above paper, is very abstract. One of the objectives of this paper is to give more manageable characterizations of the limit phase space. Under additional hypotheses on the domains $\Omega_\epsilon$ we also give a simple description of the limit equation. If, in addition, $\mathcal{M}$ is a $k$-sphere and the nonlinearity of the above equations is dissipative, then, as we will prove, for every $\epsilon > 0$ small enough the corresponding equation on $\Omega_\epsilon$ possesses an inertial manifold, i.e. an invariant manifold containing the attractor of the equation. We thus obtain the existence of inertial manifolds for reaction-diffusion equations on certain classes of thin domains of genuinely high dimension.

## 1. INTRODUCTION

In this paper we study a family of semilinear reaction-diffusion equations on spatial domains $\Omega_\epsilon$, $\epsilon > 0$, in $\mathbb{R}^\ell$ lying close to a $k$-dimensional submanifold $\mathcal{M}$ of $\mathbb{R}^\ell$. As $\epsilon \to 0^+$, the domains $\Omega_\epsilon$ shrink onto a subset of $\mathcal{M}$ in the normal direction to $\mathcal{M}$. It was proved in the previous work [14], extending earlier results from [6] and [15], that the above family has a limit equation, which is an abstract semilinear parabolic equation defined on a certain limit phase space denoted by $H_s^1(\Omega)$.

The definition of $H_s^1(\Omega)$, given in [14] is very abstract. One of the objectives of this paper is to provide more manageable characterizations of the limit phase space. Under additional hypotheses on the domains $\Omega_\epsilon$ we also give a simple description of the limit equation. If, in addition, $\mathcal{M}$ is a $k$-sphere and the nonlinearity of the above equations is dissipative, then, as we will prove, for every $\epsilon > 0$ small enough the corresponding equation on $\Omega_\epsilon$ possesses an inertial manifold, i.e. an invariant manifold containing the attractor of the equation. Thus we obtain the existence of inertial manifolds for reaction-diffusion equations on certain classes of thin domains of genuinely high dimension.

Let us now give a more detailed description of the results of this paper. Let $\ell$, $k$ and $r$ be positive integers with $r \geq 2$, $\ell \geq 2$ and $k < \ell$. Let $\mathcal{M} \subset \mathbb{R}^\ell$ be an arbitrary imbedded $k$-dimensional submanifold of $\mathbb{R}^\ell$ of class $C^r$. Note that, in the general case considered here, the manifold is *global*, i.e. $\mathcal{M}$ need not be included in a single coordinate chart. Let us also remark that we do *not* assume $\mathcal{M}$ to be orientable.

By the Tubular neighborhood theorem (cf e.g. [1]) there exists an open set $\mathcal{U}$ in $\mathbb{R}^\ell$ and a map $\phi\colon \mathcal{U} \to \mathcal{M}$ of class $C^{r-1}$ such that whenever $x \in \mathcal{U}$ and $p \in \mathcal{M}$ then $\phi(x) = p$ if and only if the vector $x - p$ is orthogonal to $T_p\mathcal{M}$; moreover, $\epsilon x + (1-\epsilon)\phi(x) \in \mathcal{U}$ for all $x \in \mathcal{U}$ and all $\epsilon \in [0,1]$.

For $\epsilon \in [0,1]$ let us define the *curved squeezing* transformation
$$\Phi_\epsilon \colon \mathcal{U} \to \mathbb{R}^\ell,$$

$$\Phi_\epsilon(x) := \epsilon x + (1-\epsilon)\phi(x) = \phi(x) + \epsilon(x - \phi(x)). \tag{1}$$

Now let $\Omega$ be an arbitrary nonempty bounded domain in $\mathbb{R}^\ell$ with Lipschitz boundary and such that $\operatorname{Cl}\Omega \subset \mathcal{U}$. For $\epsilon \in {]0,1]}$, define the *curved squeezed domain*
$$\Omega_\epsilon := \Phi_\epsilon(\Omega).$$

Let $\epsilon \in {]0,1]}$ be arbitrary, $\omega := \Omega_\epsilon$ and consider the Neumann boundary value problem

$$\begin{aligned} u_t &= \Delta u + G(u), & t > 0,\ x \in \omega, \\ \partial_\nu u &= 0, & t > 0,\ x \in \partial\omega \end{aligned} \tag{2}$$

on $\omega$. Here, $\nu$ is the exterior normal vector field on $\partial\omega$. Suppose that $G \in C^1(\mathbb{R} \to \mathbb{R})$ is dissipative in the sense that
$$\limsup_{|s| \to \infty} G(s)/s \le -\delta_0 \quad \text{for some } \delta_0 > 0.$$

Furthermore, let $G$ satisfy the growth estimate
$$|G'(s)| \le C(1 + |s|^\beta) \quad \text{for } s \in \mathbb{R},$$
where $C$ and $\beta \in [0, \infty[$ are arbitrary real constants. If $\ell > 2$, assume, in addition, that $\beta < (2^*/2) - 1$, where $2^* = 2\ell/(\ell-2)$.

This equation can be described in abstract terms as the equation
$$\dot u + \tilde A_\epsilon u = \hat G(u) \tag{3}$$
on $H^1(\Omega_\epsilon)$. Here, the operator $\tilde A_\epsilon$ is induced by the bilinear form $\tilde a_\epsilon$
$$\tilde a_\epsilon(u,v) = \int_{\Omega_\epsilon} \nabla u \cdot \nabla v \, \mathrm{d}x$$
on $H^1(\Omega_\epsilon)$ in the sense that
$$\tilde A_\epsilon u = w \text{ if and only if } \tilde a_\epsilon(u,v) = \int_{\Omega_\epsilon} uv\, \mathrm{d}x \text{ for all } v \in H^1(\Omega_\epsilon).$$

Furthermore, $\hat G(u) := G \circ u$ is the Nemitski operator defined by $G$. We can now use the change of variables $u(x) \mapsto u(\tilde x)$, where $\tilde x = \Phi_\epsilon(x)$, to transform equation (3) to the equivalent problem
$$\dot u + A_\epsilon u = \hat G(u) \tag{4}$$
on the fixed phase space $H^1(\Omega)$. Here, the operator $A_\epsilon$ is defined by the formula
$$A_\epsilon(u \circ \Phi_\epsilon) = (\tilde A_\epsilon u) \circ \Phi_\epsilon.$$

Equation (4) defines a semiflow $\pi_\epsilon$ on $H^1(\Omega)$, which possesses a global attractor $\mathcal{A}_\epsilon$.

For $x \in \mathcal{U}$ denote by $Q(x)\colon \mathbb{R}^\ell \to \mathbb{R}^\ell$ the orthogonal projection of $\mathbb{R}^\ell \cong T_p\mathbb{R}^\ell$ onto $T_p\mathcal{M}$, where $p := \phi(x)$. Let $P(x) = I - Q(x)$. Note that $P(x)$ is the orthogonal projection of $\mathbb{R}^\ell \cong T_p\mathbb{R}^\ell$ onto the orthogonal complement of $T_p\mathcal{M}$ in $T_p\mathbb{R}^\ell \cong \mathbb{R}^\ell$.

Now let us define the space

$$H^1_s(\Omega) := \{\, u \in H^1(\Omega) \mid P(x)\nabla u(x) = 0 \text{ a.e.} \,\}. \tag{5}$$

Note that $H^1_s(\Omega)$ is a closed linear subspace of the Hilbert space $H^1(\Omega)$. Let $L^2_s(\Omega)$ be the closure in $L^2(\Omega)$ of $H^1_s(\Omega)$.

It is one of the main contributions of the paper [14] to show that the family of operators $(A_\epsilon)_{\epsilon \in ]0,1]}$ converges in a strong spectral sense to a densely defined selfadjoint operator $A_0$ in $L^2_s(\Omega)$.

We can now consider the abstract parabolic equation

$$\dot{u} + A_0 u = \hat{G}(u). \tag{6}$$

on the space $H^1_s(\Omega)$, where $H^1_s(\Omega)$ is defined in (5). Equation (6) defines a semiflow $\pi_0$ on $H^1_s(\Omega)$, which possesses a global attractor $\mathcal{A}_0$.

It is proved in [14] that, as $\epsilon \to 0^+$, the linear semigroups $e^{-tA_\epsilon}$ converge in a singular sense to the semigroup $e^{-tA_0}$ and the semiflows $\pi_\epsilon$ singularly converge to $\pi_0$. Furthermore, an upper semicontinuity result is established for the family $(\mathcal{A}_\epsilon)_{\epsilon \in [0,1]}$ of attractors.

In order to precisely define the operator $A_0$ we need some notation from [14]:

define the continuous function $J_0\colon \mathcal{U} \to \mathbb{R}$ as

$$J_0(x) := |\det(D\phi(x))|_{T_{\phi(x)}\mathcal{M}}.$$

Moreover, for every $x \in \mathcal{U}$ define the linear map $S_0(x)\colon \mathbb{R}^\ell \to \mathbb{R}^\ell$ as

$$S_0(x) := \lim_{\epsilon \to 0^+} (D\Phi_\epsilon^{-1}(\Phi_\epsilon(x)) - (1/\epsilon)P(x))$$

(the limit being taken in $\mathcal{L}(\mathbb{R}^\ell, \mathbb{R}^\ell)$). It is proved in [14] that $S_0(x)$ is well-defined and the function $S_0\colon \mathcal{U} \to \mathcal{L}(\mathbb{R}^\ell, \mathbb{R}^\ell)$ is continuous.

Define the bilinear forms

$$b_0\colon L^2_s(\Omega) \times L^2_s(\Omega) \to \mathbb{R}$$

by

$$b_0(u,v) := \int_\Omega J_0(x) u(x) v(x)\, \mathrm{d}x \tag{7}$$

and

$$a_0\colon H^1_s(\Omega) \times H^1_s(\Omega) \to \mathbb{R}$$

by

$$a_0(u,v) := \int_\Omega J_0(x) \langle S_0(x)^T \nabla u(x), S_0(x)^T \nabla v(x) \rangle\, \mathrm{d}x. \tag{8}$$

$A_0$ is defined as the operator generated by the pair $(a_0, b_0)$. More precisely,

$$A_0 u = w \text{ if and only if } a_0(u,v) = b_0(w,v) \text{ for all } v \in H^1_s(\Omega).$$

Now, for $p \in \mathcal{M}$ define the *normal section* $\Omega_p$ *of* $\Omega$ *at* $p$ to be the set of all $x \in \Omega$ with $\phi(x) = p$. The first of our results (Theorem 3.1) shows that functions in $L^2_s(\Omega)$ are a.e. (relative to the corresponding Hausdorff measures) constant along the connected components of $\Omega_p$. This leads to a first characterization of the space $H^1_s(\Omega)$, see Corollary 3.2. If $\Omega$ has *connected normal sections*, i.e. if the set $\Omega_p$ is connected for all $p \in \mathcal{M}$, then Theorems 3.3 and 3.4 completely characterize the spaces $L^2_s(\Omega)$ and $H^1_s(\Omega)$. Under some additional regularity hypotheses, Theorem 3.5 and its Corollary provide a simple description of the limit operator $A_0$ and the corresponding limit equation. In particular, $A_0$ is equivalent to a relatively bounded perturbation of the Laplace-Beltrami operator on an open subset of $\mathcal{M}$. If $\mathcal{G} := \phi(\Omega)$ is a $k$-dimensional sphere and some additional hypotheses are satisfied then the eigenvalues of the limit operator $A_0$ satisfy a certain 'gap condition' (cf Theorem 3.7). We can then apply a version of the inertial manifold theorem from [16] (cf Theorem 3.8) which shows that for all small $\epsilon \geq 0$ there is an invariant manifold $\mathcal{I}_\epsilon$ for the equation

$$\dot{u} + A_\epsilon u = \hat{G}(u)$$

containing the attractor of this equation. The manifolds $\mathcal{I}_\epsilon$ (resp. the reduced equations on $\mathcal{I}_\epsilon$) converge, in a regular $C^1$-sense, to the manifold $\mathcal{I}_0$ (resp. to the reduced equations on $\mathcal{I}_0$).

The proof of Theorem 3.8 relies essentially on the 'gap condition' mentioned above. Therefore we can obtain the same result if we consider, instead of a sphere, any manifold with the property that the eigenvalues of the limit operator $A_0$ satisfy the 'gap condition'. It is worth to mention that there is a class of compact manifolds without boundary, which exibit large gaps in the spectrum of the corresponding Laplace-Beltrami operator, and for which the 'gap condition' is actually satisfied. It is the class of manifolds satisfying the following property:

all the geodesics are closed and their lenght is an integer multiple of a fixed positive number $T$.

This fact has been known for quite a long time (see e.g. [4, 8] and the references contained therein), but only recently it has been observed that such manifolds provide a new class of spatial domains on which reaction-diffusion equations possess inertial manifolds (see [10]). This is a remarkable fact, since previous results about inertial manifolds assumed that the spatial domains were segments in $\mathbb{R}$, rectangles in $\mathbb{R}^2$ or cubes in $\mathbb{R}^3$ (like in [12]), or equilateral triangles (like in [11]). Note that in [18] Temam and Wang already exploited the large gaps occurring in the spectrum of the Laplace-Beltrami operator on $\mathbb{S}^2$ (acting on 1-forms), in order to construct inertial manifolds for the Navier-Stokes equation on $\mathbb{S}^2$.

In this paper we obtain the existence of inertial manifolds for reaction-diffusion equations on certain classes of domains in $\mathbb{R}^\ell$ which are 'thin' in $\ell - k$ spatial directions but *not* thin in the remaining $k$-directions. Since we may choose $\ell$ and $k$ arbitrary with $k \leq \ell - 1$ we may therefore term these domains as being of 'genuinely high dimension'.

## 2. PRELIMINARIES

Given an arbitrary positive integer $m$, we denote by $\mathcal{H}^m$ the $m$-dimensional Hausdorff measure on $\mathbb{R}^\ell$ induced by the Euclidean metric. We need the following special case of the general coarea formula for Riemannian manifolds from [5]:

**Theorem 2.1.** *Suppose $g\colon \mathcal{U} \to \mathbb{R}$ is Lebesgue-measurable and $g \geq 0$ (resp. $g$ is Lebesgue-integrable). Then for $\mathcal{H}^k$-a.a. $p \in \mathcal{M}$ the function $g|\phi^{-1}\{p\}$ is $\mathcal{H}^{\ell-k}$-measurable (resp. $\mathcal{H}^{\ell-k}$-integrable), the function*

$$p \mapsto \int_{\phi^{-1}\{p\}} g(x)\,\mathrm{d}\mathcal{H}^{\ell-k}(x)$$

*is $\mathcal{H}^k$-measurable (resp. $\mathcal{H}^k$-integrable) and*

$$(9) \qquad \int_{\mathcal{U}} J_0(x) g(x)\,\mathrm{d}x = \int_{\mathcal{M}} \left( \int_{\phi^{-1}\{p\}} g(x)\,\mathrm{d}\mathcal{H}^{\ell-k}(x) \right) \mathrm{d}\mathcal{H}^k(p),$$

*where, as before,*

$$J_0(x) := |\det(D\phi(x))|_{T_{\phi(x)}\mathcal{M}}.$$

We will now recall a few classical definitions and results about Sobolev spaces on Riemannian manifolds. For more details on this subject, the reader is referred to [3], [7] and [17].

For the rest of this section, let $S \subset \mathcal{M}$ be open in $\mathcal{M}$. We denote by $L^2(S)$ (resp. $L^2_{\mathrm{loc}}(S)$) the set of all square integrable (resp. locally square integrable) $\mathcal{H}^k$-measurable functions defined on $S$. Besides, we denote by $\mathbb{L}^2(S)$ (resp. $\mathbb{L}^2_{\mathrm{loc}}(S)$) the space of all $\mathcal{H}^k$-measurable tangent vector fields $X$ on $S$ such that the function $p \mapsto \langle X(p), X(p) \rangle$ is integrable (resp. locally integrable) on $S$.

**Definition 2.2.** *Let $u \in C^m(S)$, $1 \leq m \leq r$. The gradient $\nabla u$ of $u$ is the $C^{m-1}$ vector field on $S$ defined by*

$$\langle \nabla u(p), h \rangle = (\,\mathrm{d}u(p), h)\quad \text{for all } p \in S \text{ and all } h \in T_p\mathcal{M}.$$

*(Here we denote by $(\cdot, \cdot)$ the duality product in $T_p\mathcal{M}$.)*

It follows that whenever $\tilde{u}\colon U \to \mathbb{R}$ is a $C^1$-extension of the function $u$ to a neighborhood $U$ of $S$ in $\mathbb{R}^\ell$ (e.g. $\tilde{u} := u \circ \phi$), then, for $p \in S$, $\nabla u(p)$ is the orthogonal projection of the usual gradient $\nabla \tilde{u}(p) \in \mathbb{R}^\ell$ onto $T_p\mathcal{M}$.

Let us denote by $\hat{\nabla}$ the Levi-Civita connection on $\mathcal{M}$. For a given vector field $X$ on $S$ of class $C^m$, $1 \leq m \leq r$, and for $p \in S$, let us define the linear map

$$R(X, p)\colon T_p\mathcal{M} \to T_p\mathcal{M}, \quad h \mapsto \hat{\nabla}_h X(p).$$

It follows again that whenever $\tilde{X}\colon U \to \mathbb{R}^\ell$ is a $C^1$-extension of the vector field $X$ to a neighborhood $U$ of $S$ in $\mathbb{R}^\ell$ (e.g. $\tilde{X} := X \circ \phi$), then, for $p \in S$, $\hat{\nabla}_h X(p)$ is the orthogonal projection of $D\tilde{X}(p)h \in \mathbb{R}^\ell$ onto $T_p\mathcal{M}$. Here, $D\tilde{X}(p)\colon \mathbb{R}^\ell \to \mathbb{R}^\ell$ is the usual Fréchet derivative of $\tilde{X}$ at $p$.

**Definition 2.3.** *Let $X$ be a vector field $X$ on $S$ of class $C^m$, $1 \leq m \leq r$. The divergence of $X$ is the $C^{m-1}$ function defined by*

$$(\operatorname{div} X)(p) := \operatorname{trace} R(X, p) \text{ for } p \in S.$$

**Definition 2.4.** *Let $u \in L^2_{\operatorname{loc}}(S)$. We say that $u \in H^1_{\operatorname{loc}}(S)$ if and only if one of the following equivalent properties is satisfied:*

1. *for every chart $\tau \colon V \subset S \to \mathcal{O} \subset \mathbb{R}^k$ the function $u \circ \tau^{-1} \colon \mathcal{O} \to \mathbb{R}$ is in $H^1_{\operatorname{loc}}(\mathcal{O})$;*
2. *there exists a sequence of functions $(\psi_n)_{n \in \mathbb{N}}$, $\psi_n \in C^1(S)$, such that, for every open set $\mathcal{V} \subset\subset S$,*

$$\int_{\mathcal{V}} |\psi_n - \psi_m|^2 \, d\mathcal{H}^k + \int_{\mathcal{V}} \langle \nabla \psi_n - \nabla \psi_m, \nabla \psi_n - \nabla \psi_m \rangle \, d\mathcal{H}^k \to 0 \quad \text{as } n, m \to \infty$$

*and such that $\psi_n \to u$ in $L^2_{\operatorname{loc}}(S)$ as $n \to \infty$;*
3. *there exists a vector field $V \in \mathbb{L}^2_{\operatorname{loc}}(S)$ such that for every vector field $\Psi$ of class $C^r$ with $\operatorname{supp} \Psi \subset\subset S$ the following equality holds:*

$$\int_S u \operatorname{div} \Psi \, d\mathcal{H}^k = -\int_S \langle V, \Psi \rangle \, d\mathcal{H}^k.$$

*We call $V$ the 'weak' gradient of $u$ and we write $V =: \nabla u$.*

*We say that $u \in H^1(S)$ if and only if $u \in H^1_{\operatorname{loc}}(S)$, $u \in L^2(S)$ and $\nabla u \in \mathbb{L}^2(S)$. For $u \in H^1(S)$, we set*

$$|u|_{H^1(S)} := \left( \int_S |u|^2 \, d\mathcal{H}^k + \int_S \langle \nabla u, \nabla u \rangle \, d\mathcal{H}^k \right)^{1/2}.$$

**Remark 2.5.** *The Cauchy type condition in (2) above implies that there exists a vector field $V \in \mathbb{L}^2_{\operatorname{loc}}(S)$ such that $\nabla \psi_n \to V$ in $\mathbb{L}^2_{\operatorname{loc}}(S)$ as $n \to \infty$. Obviously $V = \nabla u$ in the sense of property (3). It follows that $C^1$-functions are dense in $H^1_{\operatorname{loc}}(S)$.*

Let $X$ and $Y$ be two vector fields of class $C^1$ on $S$. For a given $p \in S$, let $h_j = h_j(p)$, $j = 1, \ldots, k$, be an orthonormal basis of $T_p S$, and define

$$\langle\langle X, Y \rangle\rangle(p) := \langle X(p), Y(p) \rangle + \sum_{j=1}^k \langle \hat{\nabla}_{h_j} X(p), \hat{\nabla}_{h_j} Y(p) \rangle.$$

It can be shown that $\langle\langle X, Y \rangle\rangle(p)$ does not depend on the choice of the orthonormal basis of $T_p S$, so the assignment $p \mapsto \langle\langle X, Y \rangle\rangle(p)$ defines a function on $S$, which turns out to be continuous.

**Definition 2.6.** *Let $X \in \mathbb{L}^2_{\operatorname{loc}}(S)$. We say that $X \in \mathbb{H}^1_{\operatorname{loc}}(S)$ if and only if one of the following equivalent properties is satisfied:*

1. *for every chart $\tau \colon V \subset S \to \mathcal{O} \subset \mathbb{R}^k$ and for $j = 1, \ldots, k$, the $j$-th component $X_j \colon \mathcal{O} \to \mathbb{R}$ of $X$ with respect to the coordinate system $\tau$ belongs to $H^1_{\operatorname{loc}}(\mathcal{O})$;*

2. *there exists a sequence of vector fields* $(\Psi_n)_{n\in\mathbb{N}}$, $\Psi_n \in C^1(S)$, *such that, for every open set* $\mathcal{V} \subset\subset S$,
$$\int_{\mathcal{V}} \langle\langle \Psi_n - \Psi_m, \Psi_n - \Psi_m \rangle\rangle \, \mathrm{d}\mathcal{H}^k \to 0 \quad \text{as } n, m \to \infty$$
*and such that* $\Psi_n \to X$ *in* $\mathbb{L}^2_{\mathrm{loc}}(S)$ *as* $n \to \infty$.

If $X \in \mathbb{H}^1_{\mathrm{loc}}(S)$ we can compute $\hat{\nabla}_h X$ in coordinates for almost every $p \in S$ and for all $h \in T_p \mathcal{M}$, so it makes sense to define $\langle\langle X, Y \rangle\rangle(p)$ for vector fields $X$ and $Y \in \mathbb{H}^1_{\mathrm{loc}}(S)$. We say that $X \in \mathbb{H}^1(S)$ if and only if $X \in \mathbb{H}^1_{\mathrm{loc}}(S)$ and the function $p \mapsto \langle\langle X, X \rangle\rangle(p)$ is integrable on $S$. For $X \in \mathbb{H}^1(S)$, we set
$$|X|_{\mathbb{H}^1(S)} := \left( \int_S \langle\langle X, X \rangle\rangle \, \mathrm{d}\mathcal{H}^k \right)^{1/2}.$$

**Definition 2.7.** *We say that a function* $u \in H^1(S)$ *belongs to* $H^2(S)$ *if and only if* $\nabla u \in \mathbb{H}^1(S)$.

For $u \in H^2(S)$, we define the Laplacian
$$\Delta_S u(p) := \mathrm{div}(\nabla u)(p).$$

3. The main results

In this section we will state the principal results of this paper. Most of the proofs will be given in Section 4. Let us remark that Theorems 3.1, 3.3 and Corollary 3.2 are valid without the assumption that $\Omega$ have Lipschitz boundary.

Recall that, for $p \in \mathcal{M}$, we define $\Omega_p := \{ x \in \Omega \mid \phi(x) = p \}$. Moreover, for $x \in \Omega_p$, let $\Omega_p(x)$ be the connected component of $x$ in $\Omega_p$.

We can now state the first result of this paper:

**Theorem 3.1.** *Let* $u \in L^2_{\mathrm{loc}}(\Omega)$ *and assume there exists a sequence* $(u_m)_{m\in\mathbb{N}}$ *in* $H^1_{\mathrm{loc}}(\Omega)$, *with* $P(x)\nabla u_m(x) = 0$ *a.e. in* $\Omega$ *for all* $m \in \mathbb{N}$, *such that* $u_m \to u$ *in* $L^2_{\mathrm{loc}}(\Omega)$ *as* $m \to \infty$. *Under this assumption, there exists a set* $Z \subset \mathcal{M}$, $\mathcal{H}^k(Z) = 0$, *and for all* $p \in \mathcal{M} \setminus Z$ *there exists a set* $S_p \subset \phi^{-1}(p)$, $\mathcal{H}^{\ell-k}(S_p) = 0$, *such that the following property holds:*

*for all* $p \in \mathcal{M} \setminus Z$ *and for all* $\overline{x} \in \Omega_p$ *there exists a constant* $v(p, \overline{x}) \in \mathbb{R}$ *such that* $u(x) = v(p, \overline{x})$ *for all* $x \in \Omega_p(\overline{x}) \setminus S_p$.

The assumptions above are in particular satisfied if $u \in L^2_s(\Omega)$ and a-fortiori if $u \in H^1_s(\Omega)$. Theorem 3.1 says that, up to a set of measure zero, the functions in $L^2_s(\Omega)$ are constant on each connected component of the normal section $\Omega_p$ of $\Omega$ at $p \in \mathcal{M}$.

Theorem 3.1 leads to the following simple characterization of the space $H^1_s(\Omega)$:

**Corollary 3.2.** *For* $u \in H^1(\Omega)$ *the following conditions are equivalent:*
1. $P(x)\nabla u(x) = 0$ *a.e. in* $\Omega$;
2. *There exists a set* $Z \subset \mathcal{M}$, $\mathcal{H}^k(Z) = 0$, *and for all* $p \in \mathcal{M} \setminus Z$ *there exists a set* $S_p \subset \phi^{-1}(p)$, $\mathcal{H}^{\ell-k}(S_p) = 0$, *such that the following property holds:*
   *for all* $p \in \mathcal{M} \setminus Z$ *and for all* $\overline{x} \in \Omega_p$ *there exists a constant* $v(p, \overline{x}) \in \mathbb{R}$ *such that* $u(x) = v(p, \overline{x})$ *for all* $x \in \Omega_p(\overline{x}) \setminus S_p$.

For domains $\Omega$ having connected normal sections, Theorem 3.1 implies that functions in $L_s^2(\Omega)$ depend only on the variable $p \in \mathcal{M}$. We now show that much more can be proved in this case.

Whenever $\Omega$ has connected normal sections then set $\mathcal{G} := \phi(\Omega)$ and define

$$\mu(p) := \mathcal{H}^{\ell-k}(\Omega_p) \quad \text{for } p \in \mathcal{G}.$$

The set $\mathcal{G}$ is open in $\mathcal{M}$ by the surjective mapping theorem, since $D\phi(x)\colon \mathbb{R}^\ell \to T_{\phi(x)}\mathcal{M}$ is surjective for all $x \in \mathcal{U}$. Moreover, by the coarea formula the function $\mu\colon \mathcal{G} \to \mathbb{R}$ is $\mathcal{H}^k$-measurable and, in fact, integrable on $\mathcal{G}$.

The following theorem fully characterizes the space $H_s^1(\Omega)$ when $\Omega$ has connected normal sections.

**Theorem 3.3.** *Assume that $\Omega$ has connected normal sections. Let $u \in L_s^2(\Omega)$. Then there exists an null set $S$ in $\mathbb{R}^\ell$ and a function $v \in L^2_{\mathrm{loc}}(\mathcal{G})$ such that $u(x) = v(\phi(x))$ for all $x \in \Omega \setminus S$; moreover $\mu^{1/2}v \in L^2(\mathcal{G})$. If $u \in H_s^1(\Omega)$, then $v \in H^1_{\mathrm{loc}}(\mathcal{G})$,*

$$\nabla u(x) = D\phi(x)^T \nabla v(\phi(x)) \quad \text{a.e. in } \Omega \tag{10}$$

*and $\mu^{1/2}\nabla v \in \mathbb{L}^2(\mathcal{G})$. Conversely, let $v \in L^2_{\mathrm{loc}}(\mathcal{G})$ be such that $\mu^{1/2}v \in L^2(\mathcal{G})$ and set $u(x) := v(\phi(x))$. Then $u \in L_s^2(\Omega)$. If $v \in H^1_{\mathrm{loc}}(\mathcal{G})$ and $\mu^{1/2}\nabla v \in \mathbb{L}^2(\mathcal{G})$, then $u \in H_s^1(\Omega)$.*

The main consequence of Theorem 3.3 is the following:

**Theorem 3.4.** *Suppose that $\Omega$ has connected normal sections. Define*

$$L^2(\mu, \mathcal{G}) := \{\, v \in L^2_{\mathrm{loc}}(\mathcal{G}) \mid \mu^{1/2}v \in L^2(\mathcal{G}) \,\}.$$

*Then $L^2(\mu, \mathcal{G})$, endowed with the scalar product*

$$b_\mu(v_1, v_2) := \int_\mathcal{G} \mu(p) v_1(p) v_2(p) \, d\mathcal{H}^k(p),$$

*is a Hilbert space. Moreover, define*

$$H^1(\mu, \mathcal{G}) := \{\, v \in H^1_{\mathrm{loc}}(\mathcal{G}) \mid \mu^{1/2}v \in L^2(\mathcal{G}),\ \mu^{1/2}\nabla v \in \mathbb{L}^2(\mathcal{G}) \,\}.$$

*and*

$$a_\mu(v_1, v_2) := \int_\mathcal{G} \mu(p) \langle \nabla v_1(p), \nabla v_2(p) \rangle \, d\mathcal{H}^k(p) \quad \text{for } v_1 \text{ and } v_2 \in H^1(\mu, \mathcal{G}).$$

*Then $H^1(\mu, \mathcal{G})$, endowed with the scalar product $a_\mu(\cdot, \cdot) + b_\mu(\cdot, \cdot)$, is a Hilbert space. Let $\jmath$ be the linear map*

$$\jmath\colon L_s^2(\Omega) \to L^2(\mu, \mathcal{G}), \quad u \mapsto v,$$

*where $v$ is the function given by Theorem 3.3.*

*Then $\jmath$ is an isometry of the Hilbert space $(L_s^2(\Omega), b_0(\cdot, \cdot))$ onto $L^2(\mu, \mathcal{G})$.*

*Furthermore, the restriction of the map $\jmath$ to $H_s^1(\Omega)$ is an isometry of the Hilbert space $(H_s^1(\Omega), a_0(\cdot, \cdot) + b_0(\cdot, \cdot))$ onto $H^1(\mu, \mathcal{G})$.*

*Let $A_\mu$ be the self-adjoint operator in $L^2(\mu, \mathcal{G})$ generated by the pair $(a_\mu, b_\mu)$. Then $\jmath$ restricts to an isometry $\jmath'$ of $D(A_0)$ onto $D(A_\mu)$ and $A_0 = \jmath'^{-1} A_\mu \jmath'$.*

In what follows, we denote by $\partial \mathcal{G}$ the topological boundary of $\mathcal{G}$ in $\mathcal{M}$.

**Theorem 3.5.** *Suppose that $\mathcal{G}$ is orientable (as a submanifold of $\mathcal{M}$), $\partial \mathcal{G} = \emptyset$ and the function $\mu$ is of class $C^1$ on $\mathcal{G}$.*

*Then $D(A_\mu) = H^2(\mathcal{G})$ and, for $u \in D(A_\mu)$,*

$$(A_\mu u)(p) = -(1/\mu(p))\operatorname{div}(\mu(p)\nabla u(p)) \quad \mathcal{H}^k\text{-a.e. in } \mathcal{G}.$$

*Proof.* The proof follows from the regularity theory for elliptic equations and from the divergence formula on Riemannian manifolds. Easy details are omitted. □

Theorem 3.5 clearly implies the following

**Corollary 3.6.** *Under the assumptions of Theorem 3.5, the limit equation (6) is equivalent to the following reaction-diffusion equation on $\mathcal{G}$:*

$$u_t = (1/\mu(p))\operatorname{div}(\mu(p)\nabla u) + G(u(p)), \quad t > 0, \, p \in \mathcal{G}. \quad \square$$

Instead of assuming $\partial \mathcal{G} = \emptyset$ we may alternatively assume that $\partial \mathcal{G}$ is a $k-1$-dimensional $C^2$-submanifold of $\mathcal{M}$ and that the function $\mu$ can be extended to a strictly positive $C^1$-function on $\operatorname{Cl}\mathcal{G}$. In this case it not difficult to see that the domain of the operator $A_\mu$ is the set of all functions $u \in H^2(\mathcal{G})$ satisfying the boundary condition

$$\langle \nabla u(p), \nu(p) \rangle = 0 \quad \mathcal{H}^{k-1}\text{-a.e. on } \partial \mathcal{G}$$

in the sense of traces. Here $\nu(p) \in T_p\mathcal{M}$, $p \in \partial\mathcal{G}$, is the outward normal vector field on $\partial\mathcal{G}$. Again, for $u \in D(A_\mu)$, one has

$$(A_\mu u)(p) = -(1/\mu(p))\operatorname{div}(\mu(p)\nabla u(p)) \quad \text{a.e. in } \mathcal{G}.$$

Thus the limit equation (6) takes the form

$$u_t = (1/\mu(p))\operatorname{div}(\mu(p)\nabla u) + G(u(p)), \qquad t > 0, \, p \in \mathcal{G},$$
$$\langle \nabla u(p), \nu(p) \rangle = 0, \qquad t > 0, \, p \in \partial\mathcal{G}.$$

We will now see that, for thin domains close to spheres, a spectral gap condition is satisfied, which can be used to prove existence of inertial manifolds:

**Theorem 3.7.** *Suppose $\Omega$ has connected normal sections, regard $\mathbb{R}^{k+1}$ as isometrically imbedded into $\mathbb{R}^\ell$, let $r \in \,]0,\infty[$ be arbitrary and assume that*

$$\mathcal{G} = \mathbb{S}^k(r) := \{\, x \in \mathbb{R}^{k+1} \mid \langle x, x \rangle = r^2 \,\}$$

*(i.e. $\mathcal{G}$ the $k$-dimensional sphere in $\mathbb{R}^\ell$ of radius $r$ centered at $0$). Suppose that*

$$C_\mu := \sup_{p \in \mathbb{S}^k(r)} (1/\mu(p))\langle \nabla\mu(p), \nabla\mu(p) \rangle^{1/2} \leq 1/(4r)^2.$$

*Under these assumptions the repeated sequence $(\lambda_j^0)_{j \in \mathbb{N}}$ of the eigenvalues of the limit operator $A_0$ satisfies the following 'gap' condition:*

$$(11) \qquad \limsup_{\nu \to \infty} \frac{\lambda_{\nu+1}^0 - \lambda_\nu^0}{(\lambda_\nu^0)^{1/2}} > 0.$$

We will now state an inertial manifold theorem established in the paper [16]. To this end, we need some notation.

For every $\epsilon \in [0,1]$ denote by $(\lambda_j^\epsilon)_{j\in\mathbb{N}}$ the repeated sequence of eigenvalues of the operator $A_\epsilon$ and by $(w_j^\epsilon)_{j\in\mathbb{N}}$ a corresponding orthonormal sequence of eigenvalues.

For every $\nu \in \mathbb{N}$ let $X_{\epsilon,\nu,1}$ be the span of the vectors $w_j^\epsilon$, $j = 1, \ldots, \nu$ and let $X_{\epsilon,\nu,2}$ be the orthogonal complement of $X_{\epsilon,\nu,1}$ in $L^2(\Omega)$ if $\epsilon > 0$ and in $L_s^2(\Omega)$ if $\epsilon = 0$. Let $A_{\epsilon,\nu,i}$ be the restriction of $A_\epsilon$ to $X_{\epsilon,\nu,i}$ for $i = 1, 2$. Let $E_{\epsilon,\nu}\xi := \sum_{j=1}^\nu \xi_j w_j^\epsilon$, $\xi \in \mathbb{R}^\nu$ and $P_{\epsilon,\nu,i}$ be the orthogonal projection of $L^2(\Omega)$ onto $X_{\epsilon,\nu,i}$, $i = 1, 2$ if $\epsilon > 0$ and $P_{\epsilon,\nu,i}$ be the orthogonal projection of $L_s^2(\Omega)$ onto $X_{\epsilon,\nu,i}$, $i = 1, 2$ if $\epsilon = 0$.

Finally, whenever $\epsilon \in [0,1]$ and $F\colon H^1(\Omega) \to L^2(\Omega)$ is a locally Lipschitzian (nonlinear) operator mapping $H_s^1(\Omega)$ into $L^2(\Omega)$, then $\pi_{\epsilon,F}$ is the local semiflow on $H^1(\Omega)$ for $\epsilon > 0$ and on $H_s^1(\Omega)$ for $\epsilon = 0$ generated by the solutions of the equation

$$\dot{u} + A_\epsilon u = F(u).$$

**Theorem 3.8.** *Suppose the eigenvalues of $A_0$ satisfy the following gap condition:*

$$\text{(12)} \qquad \limsup_{\nu \to \infty} \frac{\lambda_{\nu+1}^0 - \lambda_\nu^0}{(\lambda_\nu^0)^{1/2}} > 0.$$

*Then there are an $\epsilon_0 > 0$ and an open bounded set $U \subset H^1(\Omega)$ such that for every $\epsilon \in [0, \epsilon_0[$ the attractor $\mathcal{A}_\epsilon$ of the semiflow $\pi_{\epsilon,\hat{f}}$ lies in $U$.*

*Furthermore, there exists a globally Lipschitzian map $g \in C^1(H^1(\Omega) \to L^2(\Omega))$ with $g(u) = \hat{f}(u)$ for $u \in U$.*

*Besides, there is a positive integer $\nu$ and for every $\epsilon \in [0, \epsilon_0[$ there is a map $\Lambda_\epsilon \in C^1(\mathbb{R}^\nu \to H^1(\Omega))$ if $\epsilon > 0$ and $\Lambda_\epsilon \in C^1(\mathbb{R}^\nu \to H_s^1(\Omega))$ if $\epsilon = 0$ such that*

$$\text{(13)} \qquad P_{\epsilon,\nu,1} \circ \Lambda_\epsilon = E_{\epsilon,\nu}$$

*and $\mathcal{I}_\epsilon := \Lambda_\epsilon(\mathbb{R}^\nu)$ is a $C^1$-manifold which is invariant with respect to the semiflow $\pi_{\epsilon,g}$.*

*Finally, there is an open set $V \subset \mathbb{R}^\nu$ such that for every $\epsilon \in [0, \epsilon_0[$*

$$\mathcal{A}_\epsilon \subset \Lambda_\epsilon(V) \subset U$$

*and the set $\Lambda_\epsilon(V)$ is positively invariant with respect to the semiflow $\pi_{\epsilon,\hat{f}}$.*

*The reduced equation on $\Lambda_\epsilon(\mathbb{R}^\nu)$ takes the form*

$$\text{(14)} \qquad \dot{\xi} = v_\epsilon(\xi), \quad \xi \in \mathbb{R}^\nu,$$

*where*

$$v_\epsilon \colon \mathbb{R}^\nu \to \mathbb{R}^\nu, \quad \xi \mapsto -A_\epsilon E_{\epsilon,\nu}\xi + P_{\epsilon,\nu,1} g(\Lambda_\epsilon(\xi)).$$

*Moreover, whenever $\epsilon_n \to 0^+$ and $\xi_n \to \xi_0$ in $\mathbb{R}^\nu$, then*

$$\text{(15)} \qquad |\Lambda_{\epsilon_n}(\xi_n) - \Lambda_0(\xi_0)|_{\epsilon_n} + \sum_{j=1}^\nu |\partial_j \Lambda_{\epsilon_n}(\xi_n) - \partial_j \Lambda_0(\xi_0)|_{\epsilon_n} \to 0$$

*and*

$$\text{(16)} \qquad |v_{\epsilon_n}(\xi_n) - v_0(\xi_0)|_{\mathbb{R}^\nu} + \sum_{j=1}^\nu |\partial_j v_{\epsilon_n}(\xi_n) - \partial_j v_0(\xi_0)|_{\mathbb{R}^\nu} \to 0.$$

□

**Remark 3.9.** *Our inertial manifold $\mathcal{I}_\epsilon$ is (globally) invariant with respect to the modified semiflow $\pi_{\epsilon,g}$, which coincides with the original semiflow $\pi_{\epsilon,\hat{f}}$ on the neighborhood $U$ of the attractor $\mathcal{A}_\epsilon$. Thus, close to the attractor, $\mathcal{I}_\epsilon$ is a locally invariant manifold for the 'true' semiflow $\pi_{\epsilon,\hat{f}}$. In comparison with similar results contained in the literature (e.g. in [6]), Theorem 3.8 seems to be sharper. Infact, in order to prove existence of inertial manifolds, one usually finds $L^\infty$-estimates for the attractors and then modifies the nonlinearity $f\colon \mathbb{R} \to \mathbb{R}$ so as to obtain a bounded nonlinearity $\tilde{f}\colon \mathbb{R} \to \mathbb{R}$ which induces a globally Lipschitzian Nemitski operator from $H^1$ to itself. As a consequence of this modification of the function $f$ (instead of the Nemitski operator $\hat{f}$), one has that the modified semiflow coincides with the original one only on the attractor, while it is different from the latter on everey neighborhood of the attractor.*

**Remark 3.10.** *Actually Theorem 3.8 was stated and proved in [16] for the 'flat' squeezing case considered there. However, the largely abstract proof given in [16] carries over almost verbatim to the present more general situation. Trivial modifications are left to the reader.*

Combining Theorems 3.7 and 3.8 we arrive at the following important

**Corollary 3.11.** *Assume the hypotheses of Theorem 3.7. Then the conclusions of Theorem 3.8 hold.* □

We thus obtain the existence of inertial manifolds for reaction-diffusion equations on certain classes of domains in $\mathbb{R}^\ell$ which are 'thin' in $\ell - k$ spatial directions but *not* thin in the remaining $k$-directions. Since we may choose $\ell$ and $k$ arbitrary with $k \leq \ell - 1$ we may therefore term these domains as being of 'genuinely high dimension'.

## 4. The proofs

For every $p \in \mathcal{M}$ there is an open set $V_p$ in $\mathcal{M}$, $p \in V_p$, a chart $\tau = \tau_p \colon V_p \to \mathbb{R}^k$ of $\mathcal{M}$ and $C^{r-1}$-maps $\nu_j = \nu_{p,j}\colon V_p \to \mathbb{R}^\ell$, $j = 1, \ldots, \ell - k$ such that for every $q \in V_p$ the vectors $\nu_j(q)$, $j = 1, \ldots, \ell - k$, form an orthonormal basis of the orthogonal complement of $T_q\mathcal{M}$ in $T_q\mathbb{R}^\ell \cong \mathbb{R}^\ell$. For $j = 1, \ldots, \ell - k$ define the function $\alpha_j = \alpha_{p,j}\colon \phi^{-1}(V_p) \to \mathbb{R}$ by

$$\alpha_j(x) = \langle x - \phi(x), \nu_j(x)\rangle, \text{ for } x \in \phi^{-1}(V_p).$$

Set $\alpha(x) := (\alpha_1(x), \ldots, \alpha_{\ell-k}(x))$ for $x \in \phi^{-1}(V_p)$. It is easily proved that the map $\Gamma = \Gamma_p \colon \phi^{-1}(V_p) \to \mathbb{R}^k \times \mathbb{R}^{\ell-k}$, $x \mapsto (\xi, s)$ where $\xi = \tau(\phi(x))$ and $s = \alpha(x)$, is a $C^{r-1}$-diffeomorphism of $\phi^{-1}(V_p)$ onto an open set $\mathcal{O} = \mathcal{O}_p$ in $\mathbb{R}^k \times \mathbb{R}^{\ell-k} = \mathbb{R}^\ell$. The inverse map $\zeta \colon \mathcal{O} \to \phi^{-1}(V_p)$ is given by

$$\zeta\colon (\xi, s) \mapsto \sigma(\xi) + \sum_{j=1}^{\ell-k} s_j \nu_j(\sigma(\xi)) \text{ for } (\xi, s) \in \mathcal{O}.$$

Here, $\sigma := \tau^{-1}$.

**Proposition 4.1.** *Let $E \subset \mathcal{O}_p$ be open and $F := \zeta(E)$.*

*Then $u \in H^1_{\mathrm{loc}}(F)$ if and only if $\tilde{u} := u \circ \zeta \in H^1_{\mathrm{loc}}(E)$.*

*In this case the following properties are equivalent:*

1. $P(x)\nabla u(x) = 0$ *a.e. in $E$;*
2. *For every $i = 1, \ldots, \ell - k$, $\partial_{s_i}\tilde{u}(\xi, s) = 0$ a.e. in $F$.*

*Proof.* The proof is obtained by using the well-known change of variable formula in Sobolev spaces (cf Proposition IX.6 in [2]). In fact this result implies that the first part of the proposition is true and that the following chain rule holds:

$$(17) \quad \frac{\partial \tilde{u}}{\partial s_i}(\xi, s) = \sum_{l=1}^{\ell} \frac{\partial u}{\partial x_l}(\zeta(\xi, s)) \frac{\partial \zeta_l}{\partial s_i}(\xi, s) \text{ for a.a. } (\xi, s) \in E.$$

Since

$$\frac{\partial \zeta_l}{\partial s_i}(\xi, s) = (\nu_i(\sigma(\xi)))_l$$

for $i = 1, \ldots, \ell - k$ and $l = 1, \ldots, \ell$, it follows from (17) that

$$\frac{\partial \tilde{u}}{\partial s_i}(\xi, s) = \langle \nabla u(\zeta(\xi, s)), \nu_i(\sigma(\xi)) \rangle = \langle P(\zeta(\xi, s))\nabla u(\zeta(\xi, s)), \nu_i(\sigma(\xi)) \rangle.$$

This completes the proof. □

**Proposition 4.2.** *Assume that $E \subset \mathcal{O}_p$ has the special form $E := E_1 \times E_2$ where*

$$(18) \quad E_1 := \prod_{i=1}^{k} ]a_i, b_i[ \subset \mathbb{R}^k,$$

$$(19) \quad E_2 := \prod_{j=1}^{\ell-k} ]c_j, d_j[ \subset \mathbb{R}^{\ell-k}.$$

*and set $F := \zeta(E)$ as before. Let $u \in H^1_{\mathrm{loc}}(F)$ be such that $P(x)\nabla u(x) = 0$ a.e. in $F$. Then there exist a null set $S \subset \mathbb{R}^\ell$ and a function $\tilde{v} \in H^1_{\mathrm{loc}}(E_1)$ such that $u(x) = (\tilde{v} \circ \tau \circ \phi)(x)$ for all $x \in F \setminus S$. Equivalently, set $\mathcal{E}_1 := \sigma(E_1)$ and $v(q) := (\tilde{v} \circ \tau)(q)$ for $q \in \mathcal{E}_1$. Then $v \in H^1_{\mathrm{loc}}(\mathcal{E}_1)$ and $u(x) = (v \circ \phi)(x)$ for all $x \in F \setminus S$.*

*Proof.* Let $\tilde{u} := u \circ \phi$. By Proposition 4.1

$$\frac{\partial \tilde{u}}{\partial s_i}(\xi, s) = 0 \quad \text{a.e } E, \ i = 1, \ldots, \ell - k.$$

By Lemma 2.3 in [15] it follows that there exist a null set $\tilde{S} \subset \mathbb{R}^\ell$ and a function $\tilde{v} \in H^1_{\mathrm{loc}}(E_1)$ such that

$$\tilde{u}(\xi, s) = \tilde{v}(\xi) \quad \text{for all } (\xi, s) \in E \setminus \tilde{S}.$$

Set $S := \zeta(\tilde{S} \cap E)$. Then for $x \in F \setminus S$, we have

$$u(x) = \tilde{u}(\Gamma(x)) = \tilde{v}(\tau(\phi(x))).$$

The second part of the proposition follows from Definition 2.4. □

**Proposition 4.3.** *Let $E$ and $F$ be as in Proposition 4.2. Let $u \in L^2_{\text{loc}}(F)$ and assume there exists a sequence $(u_m)_{m \in \mathbb{N}}$ in $H^1_{\text{loc}}(F)$ with $P(x)\nabla u_m(x) = 0$ a.e. in $F$ for all $m \in \mathbb{N}$, such that $u_m \to u$ in $L^2_{\text{loc}}(F)$ as $m \to \infty$. Then there exist a null set $S \subset \mathbb{R}^\ell$ and a function $\tilde{v} \in L^2_{\text{loc}}(E_1)$ such that $u(x) = (\tilde{v} \circ \tau \circ \phi)(x)$ for all $x \in F \setminus S$. Equivalently, set $\mathcal{E}_1 := \sigma(E_1)$ and $v(q) := (\tilde{v} \circ \tau)(q)$ for $q \in \mathcal{E}_1$. Then $v \in L^2_{\text{loc}}(\mathcal{E}_1)$ and $u(x) = v \circ \phi(x)$ for all $x \in F \setminus S$.*

*Proof.* By Proposition 4.1
$$\tilde{u}_m := u_m \circ \zeta \in H^1_{\text{loc}}(E).$$
Moreover, since $\zeta$ is a diffeomorphism,
$$\tilde{u} := u \circ \zeta \in L^2_{\text{loc}}(E)$$
and
$$\tilde{u}_m \to \tilde{u} \quad \text{in } L^2_{\text{loc}}(E) \text{ as } m \to \infty.$$
By Proposition 4.1 we obtain, for all $m \in \mathbb{N}$ and all $i = 1, \ldots, \ell - k$,
$$\frac{\partial \tilde{u}_m}{\partial s_i}(\xi, s) = 0 \quad \text{a.e. in } E.$$
It follows that, for all $i = 1, \ldots, \ell - k$, $\partial_{s_i}\tilde{u}(\xi, s) = 0$ in the distributional sense. Thus Lemma 2.3 in [16] again implies that there exist a null set $\tilde{S} \subset \mathbb{R}^\ell$ and a function $\tilde{v} \in L^2_{\text{loc}}(E_1)$ such that
$$\tilde{u}(\xi, s) = \tilde{v}(\xi) \quad \text{for all } (\xi, s) \in E \setminus \tilde{S}.$$
Set $S := \zeta(\tilde{S} \cap E)$. Then for $x \in F \setminus S$, we have
$$u(x) = \tilde{u}(\Gamma(x)) = \tilde{v}(\tau(\phi(x))).$$
The second part of the proposition is obvious. □

Before stating the next result, let us notice that, whenever $S$ is an $\mathcal{H}^k$-measurable subset of $\mathcal{M}$ then $\phi^{-1}(S)$ is Lebesgue measurable in $\mathbb{R}^\ell$. This follows from the fact that the map $\phi \colon \mathcal{U} \to \mathcal{M}$ is a submersion, so in local charts it can be described as the canonical projection $\pi$ of $\mathbb{R}^\ell = \mathbb{R}^k \times \mathbb{R}^{\ell-k}$ onto $\mathbb{R}^k$. Thus the above statement boils down to proving that whenever $A$ is Lebesgue measurable in $\mathbb{R}^k$ then $\pi^{-1}(A)$ is Lebesgue measurable in $\mathbb{R}^\ell$. However, this latter statement is well known to be true. In particular, if $v$ is a measurable function defined on $S$ then $u = v \circ \phi$ is a measurable function defined on $\phi^{-1}(S)$. We will use this remark implicitly in the proofs to follow.

**Proposition 4.4.** *Let $V$ be open in $\mathcal{M}$ and $U := \phi^{-1}(V)$. Suppose $v \in H^1_{\text{loc}}(V)$ and $u = v \circ \phi$ a.e. in $U$. Then $u \in H^1_{\text{loc}}(U)$ and*
$$\nabla u(x) = D\phi(x)^T \nabla v(\phi(x)) \text{ a.e. in } U.$$

*Proof.* Let us first assume that $v \in C^1(V)$. In this case we can assume w. l. o. g. that $u = v \circ \phi \in C^1(U)$. Let $x \in U$ and $h \in \mathbb{R}^\ell$ be arbitrary and set $q := \phi(x)$. Now $\phi \circ \phi = \phi$ so $u \circ \phi = u$. Since $D\phi(x)h \in T_q\mathcal{M}$ and $Q(q)$ is the orthogonal projection of $\mathbb{R}^\ell$ onto $T_q\mathcal{M}$ it follows that

$$\langle \nabla u(x), h \rangle = Du(x)h = D(u \circ \phi)(x)h = Du(q)D\phi(x)h = \langle \nabla u(q), D\phi(x)h \rangle$$
$$= \langle Q(q)\nabla u(q), D\phi(x)h \rangle = \langle \nabla v(q), D\phi(x)h \rangle = \langle D\phi(x)^T \nabla v(q), h \rangle.$$

Since $h$ is arbitrary, we see that $\nabla u(x) = D\phi(x)^T \nabla v(\phi(x))$ for all $x \in U$.

Assume now that $v \in H^1_{\mathrm{loc}}(V)$. Take a sequence $(v_m)_{m \in \mathbb{N}}$ in $C^1(V)$ such that $v_m \to v$ in $H^1_{\mathrm{loc}}(V)$ as $m \to \infty$. Let $U^0$ be an arbitrary open set with $U^0 \subset\subset U$ and let $V^0 = \phi(U^0)$. Notice that $V^0 \subset\subset V$. By (9), we have

$$\int_{U^0} J_0(x)|v_m(\phi(x)) - v(\phi(x))|^2 \, \mathrm{d}x$$
$$= \int_{V^0} \mathcal{H}^{\ell-k}(\phi^{-1}\{p\} \cap U^0)|v_m(p) - v(p)|^2 \, \mathrm{d}\mathcal{H}^k(p) \to 0$$
$$\text{as } m \to \infty.$$

Since $\inf_{x \in U^0} J_0(x) > 0$ we infer that

$$v_m \circ \phi \to v \circ \phi \quad \text{in } L^2(U^0) \text{ as } m \to \infty.$$

In the same way we have

$$\int_{U^0} J_0(x)\langle \nabla v_m(\phi(x)) - \nabla v(\phi(x)), \nabla v_m(\phi(x)) - \nabla v(\phi(x))\rangle \, \mathrm{d}x$$
$$= \int_{V^0} \mathcal{H}^{n-k}(\phi^{-1}\{p\} \cap U^0)\langle \nabla v_m(p) - \nabla v(p), \nabla v_m(p) - \nabla v(p)\rangle \, \mathrm{d}\mathcal{H}^k(p) \to 0$$
$$\text{as } m \to \infty.$$

Again we get

$$\nabla v_m \circ \phi \to \nabla v \circ \phi \quad \text{in } L^2(U^0, \mathbb{R}^\ell) \text{ as } m \to \infty.$$

Now choose an arbitrary function $\psi \in C_0^\infty(U)$ with $\operatorname{supp} \psi \subset U^0$. From what we have proved thus far, we obtain

$$\int_U u(x)\nabla\psi(x) \, \mathrm{d}x = \int_U v(\phi(x))\nabla\psi(x) \, \mathrm{d}x = \int_{U^0} v(\phi(x))\nabla\psi(x) \, \mathrm{d}x$$
$$= \lim_{m \to \infty} \int_{U^0} v_m(\phi(x))\nabla\psi(x) \, \mathrm{d}x = -\lim_{m \to \infty} \int_{U^0} \psi(x)D\phi(x)^T \nabla v_m(\phi(x)) \, \mathrm{d}x$$
$$= -\int_{U^0} \psi(x)D\phi(x)^T \nabla v(\phi(x)) \, \mathrm{d}x = -\int_U \psi(x)D\phi(x)^T \nabla v(\phi(x)) \, \mathrm{d}x$$

Since $U^0$ is arbitrary, we thus conclude both that $u \in H^1_{\mathrm{loc}}(U)$ and that $\nabla u(x) = D\phi(x)^T \nabla v(\phi(x))$ a.e. in $U$. $\square$

*Proof of Theorem 3.1.* We follow the proof of Theorem 2.5 in [15]. There exists a sequence $(p_i)_{i \in \mathbb{N}}$ of (not necessarily pairwise distinct) points in $\mathcal{M}$, a sequence $\tau_i \colon V_i = V_{p_i} \to \mathbb{R}^k$, $i \in \mathbb{N}$, of charts with inverses $\sigma_i$, a sequence $\zeta_i \colon \mathcal{O}_i = \mathcal{O}_{p_i} \subset \mathbb{R}^k \times \mathbb{R}^{\ell-k} \to \phi^{-1}(V_i)$, $i \in \mathbb{N}$, of diffeomorphisms and a sequence of sets $(F_i)_{i \in \mathbb{N}}$ such that:

1. $\Omega = \cup_{i \in \mathbb{N}} F_i$;

2. for all $i \in \mathbb{N}$, $F_i$ has the form of the set $F$ of Proposition 4.2, i.e. $F_i = \zeta_i(E_{1,i} \times E_{2,i})$, where $E_{1,i}$ and $E_{2,i}$ are products of open intervals, $i \in \mathbb{N}$.

Set $\mathcal{E}_{1,i} := \sigma_i(E_{1,i})$, $i \in \mathbb{N}$. Then, for all $i$ and $m \in N$, we have that

$$P(x)\nabla u_m(x) = 0 \quad \text{a.e. in } F_i$$

and

$$u_m|_{F_i} \to u|_{F_i} \quad \text{in } L^2_{\text{loc}}(F_i) \text{ as } m \to \infty.$$

By Proposition 4.3, for all $i \in \mathbb{N}$ there exist a null set $S_i \subset \mathbb{R}^\ell$, and a function $v_i \in L^2_{\text{loc}}(\mathcal{E}_{1,i})$ such that

$$u(x) = v_i(\phi(x)) \quad \text{for all } x \in F_i \setminus S_i.$$

Set

$$S := (\cup_{i \in \mathbb{N}} S_i) \cap \mathcal{U} \quad \text{and} \quad S_p := S \cap \phi^{-1}\{p\} \quad \text{for all } p \in \mathcal{M}.$$

Let $\chi_S$ denote the characteristic function of the set $S$. Then, by the coarea formula (9), we have

$$0 = \int_\mathcal{U} J_0(x)\chi_S(x)\,\mathrm{d}x = \int_\mathcal{M} \mathcal{H}^{\ell-k}(S_p)\,\mathrm{d}\mathcal{H}^k(p).$$

It follows that there exists a set $Z \subset \mathcal{M}$, $\mathcal{H}^k(Z) = 0$, such that, for all $p \in \mathcal{M} \setminus Z$, $\mathcal{H}^{\ell-k}(S_p) = 0$. Let $p \in \mathcal{M} \setminus Z$ be arbitrary. and fix $\overline{x} \in \Omega_p$. There exists an $i \in \mathbb{N}$ such that $\overline{x} \in F_i$. Set

$$A := \{\, z \in \Omega_p(\overline{x}) \mid \text{there exists a set } V_z \subset \Omega_p(\overline{x}),$$
$$V_z \text{ open in } \Omega_p(\overline{x}), z \in V_z, \text{ such that } u(x) = v_i(\phi(\overline{x})) \text{ for all } x \in V_z \setminus S_p \,\}.$$

Then, by Proposition 4.3, $\overline{x} \in A$, with $V_{\overline{x}} = \zeta_i(\{\tau_i(p)\} \times E_{2,i})$. If $z \in A$, then obviously $V_z \subset A$, so $A$ is open in $\Omega_p(\overline{x})$. If $z \in \Omega_p(\overline{x}) \setminus A$, then $z \in F_j$ for some $j \in \mathbb{N}$. By Proposition 4.3, it follows that

$$u(x) = v_j(\phi(\overline{x})) \quad \text{for all } x \in \zeta_j(\{\tau_j(p)\} \times E_{2,j}) \setminus S_p.$$

If $v_j(\phi(\overline{x})) = v_i(\phi(x))$, then $z \in A$, a contradiction. Thus necessarily $v_j(\phi(\overline{x})) \neq v_i(\phi(x))$ and so

$$A \cap \zeta_j(\{\sigma_j^{-1}(p)\} \times E_{2,j}) = \emptyset.$$

This implies that $A$ is closed in $\Omega_p(\overline{x})$. Since $\Omega_p(\overline{x})$ is connected, it follows that $A = \Omega_p(\overline{x})$. Therefore $\Omega_p(\overline{x}) = \cup_{z \in A} V_z$. Let $x \in \Omega_p(\overline{x}) \setminus S_p$. Then $x \in V_z$ for some $z \in A$ and hence $u(x) = v_i(\phi(\overline{x}))$. We conclude that $u(x) = v_i(\phi(\overline{x}))$ for all $x \in \Omega_p(\overline{x}) \setminus S_p$. □

*Proof of Corollary 3.2.* Condition (1) implies condition (2) by Theorem 3.1. Now suppose that condition (2) holds. Let $(F_i)_{i \in \mathbb{N}}$ and $(\zeta_i)_{i \in \mathbb{N}}$ be as in the proof of Theorem 3.1. For every $i \in \mathbb{N}$ set $\tilde{u} = u \circ \zeta_i$. Condition (2) implies that there is a null set $\tilde{Z}$ in $\mathbb{R}^k$ and for every $\xi \in E_{1,j} \setminus \tilde{Z}$ there is a null set $\tilde{S}_\xi \in \mathbb{R}^{\ell-k}$ with the property that for every $\xi \in E_{1,j} \setminus \tilde{Z}$ there is a constant $\tilde{v}(\xi)$ such that $\tilde{u}(\xi, s) = \tilde{v}(\xi)$ for all $s \in E_{2,j} \setminus \tilde{S}_\xi$. Therefore Theorem 2.5 in [15] implies that $\partial_{s_j}\tilde{u} = 0$ a.e. in $E_{1,i} \times E_{2,i}$, $j = 1, \ldots, \ell - k$. Thus Proposition 4.1 shows that $P(x)\nabla u(x) = 0$ a.e. in $F_i$, for all $i \in \mathbb{N}$. Hence $P(x)\nabla u(x) = 0$ a.e. in $\Omega$. The corollary is proved. □

*Proof of Theorem 3.3.* Let $(F_i)_{i\in\mathbb{N}}$ be as in the proof of Theorem 3.1. As in that proof, for every $i \in \mathbb{N}$ there exist a null set $S_i \subset \mathbb{R}^\ell$ and a function $v_i \in L^2_{\text{loc}}(\mathcal{E}_{1,i})$ ($v_i \in H^1_{\text{loc}}(\mathcal{E}_{1,i})$ in case $u \in H^1_s(\Omega)$) such that

$$u(x) = v_i(\phi(x)) \quad \text{for all } x \in F_i \setminus S_i.$$

Observe that $\mathcal{G} = \cup_{i\in\mathbb{N}}\mathcal{E}_{1,i}$. Again, set

$$S := (\cup_{i\in\mathbb{N}}S_i) \cap \mathcal{U} \quad \text{and} \quad S_p := S \cap \phi^{-1}\{p\} \quad \text{for all } p \in \mathcal{M}.$$

It follows that there exists a set $Z \subset \mathcal{M}$, $\mathcal{H}^k(Z) = 0$, such that, for all $p \in \mathcal{M} \setminus Z$, $\mathcal{H}^{\ell-k}(S_p) = 0$.

Let $p \in \mathcal{G} \setminus Z$ and take any $x \in \Omega_p \setminus S_p$. Then

$$v_i(p) = v_i(\phi(x)) = u(x) \text{ for all } i \in \mathbb{N} \text{ with } p \in \mathcal{E}_{1,i}.$$

In particular we obtain that

$$v_i(p) = v_j(p) \quad \text{for all } p \in (\mathcal{E}_{1i} \cap \mathcal{E}_{1j}) \setminus Z \text{ and for all } i, j \in \mathbb{N}.$$

Let us define the function $v \colon \mathcal{G} \setminus Z \to \mathbb{R}$ by

$$v(p) := v_i(p) \quad \text{if } p \in \mathcal{E}_{1i}.$$

Then $v$ is defined unambiguously and we can extend it trivially to the whole of $\mathcal{G}$. Obviously $v|_{\mathcal{E}_{1,i}} = v_i$ a.e. in $\mathcal{E}_{1,i}$, so $v \in L^2_{\text{loc}}(\mathcal{G})$. Moreover, $u(x) = v(\phi(x))$ for all $x \in \Omega \setminus S$. If $u \in H^1_s(\Omega)$ then, by Propositions 4.2 and 4.4, we have $v \in H^1_{\text{loc}}(\mathcal{G})$ and $\nabla u(x) = D\phi(x)^T \nabla v(\phi(x))$ a.e. in $\Omega$. Finally, by the coarea formula, we have:

$$\int_{\mathcal{G}} \mu(p)|v(p)|^2 \, d\mathcal{H}^k(p) = \int_{\phi(\Omega)} \mathcal{H}^{\ell-k}(\phi^{-1}(p) \cap \Omega)|v(p)|^2 \, d\mathcal{H}^k(p)$$

$$= \int_\Omega J_0(x)|v(\phi(x))|^2 \, dx = \int_\Omega J_0(x)|u(x)|^2 \, dx < \infty.$$

It is proved in [14] that $S_0(x)^T D\phi(x)^T h = h$ for all $x \in \mathcal{U}$ and all $h \in T_{\phi(x)}\mathcal{M}$. Using this, we further obtain

$$\int_{\mathcal{G}} \mu(p)\langle \nabla v(p), \nabla v(p)\rangle \, d\mathcal{H}^k(p)$$

$$= \int_{\phi(\Omega)} \mathcal{H}^{\ell-k}(\phi^{-1}(p) \cap \Omega)\langle \nabla v(p), \nabla v(p)\rangle \, d\mathcal{H}^k(p) = \int_\Omega J_0(x)\langle \nabla v(\phi(x)), \nabla v(\phi(x))\rangle \, dx$$

$$= \int_\Omega J_0(x)\langle S_0(x)^T D\phi(x)^T \nabla v(\phi(x)), S_0(x)^T D\phi(x)^T \nabla v(\phi(x))\rangle \, dx$$

$$= \int_\Omega J_0(x)\langle S_0(x)^T \nabla u(x), S_0(x)^T \nabla u(x)\rangle \, dx < \infty.$$

This completes the proof of the first part of the theorem.

Assume now that $v \in L^2_{\text{loc}}(\mathcal{G})$ and $\mu^{1/2} v \in L^2(\mathcal{G})$. Set $u(x) := v(\phi(x))$ for $x \in \Omega$. By the coarea formula we have

$$\int_\Omega J_0(x)|u(x)|^2 \, \mathrm{d}x = \int_\Omega J_0(x)|v(\phi(x))|^2 \, \mathrm{d}x$$

$$= \int_{\phi(\Omega)} \mathcal{H}^{\ell-k}(\phi^{-1}(p) \cap \Omega)|v(p)|^2 \, \mathrm{d}\mathcal{H}^k(p) = \int_\mathcal{G} \mu(p)|v(p)|^2 \, \mathrm{d}\mathcal{H}^k(p) < \infty.$$

Since $\inf_{x \in \Omega} J_0(x) > 0$, this implies that $u \in L^2(\Omega)$. Set

$$C := \sup_{x \in \Omega} \|D\phi(x)^T\|_{\mathcal{L}(\mathbb{R}^\ell, \mathbb{R}^\ell)} < \infty.$$

If $v \in H^1_{\text{loc}}(\mathcal{G})$ and $\mu^{1/2} \nabla v \in \mathbb{L}^2(\mathcal{G})$, then we obtain

$$\int_\Omega J_0(x) \langle \nabla u(x), \nabla u(x) \rangle \, \mathrm{d}x$$

$$= \int_\Omega J_0(x) \langle D\phi(x)^T \nabla v(\phi(x)), D\phi(x)^T \nabla v(\phi(x)) \rangle \, \mathrm{d}x$$

$$\leq C \int_\Omega J_0(x) \langle \nabla v(\phi(x)), \nabla v(\phi(x)) \rangle \, \mathrm{d}x = C \int_{\phi(\Omega)} \mathcal{H}^{\ell-k}(\phi^{-1}(p) \cap \Omega) \langle \nabla v(p), \nabla v(p) \rangle \, \mathrm{d}\mathcal{H}^k(p)$$

$$= C \int_\mathcal{G} \mu(p) \langle \nabla v(p), \nabla v(p) \rangle \, \mathrm{d}\mathcal{H}^k(p) < \infty.$$

Thus again $u \in H^1(\Omega)$. Since

$$\nabla u(x) = D\phi(x)^T \nabla v(\phi(x)) \quad \text{a.e. in } \Omega,$$

it follows that $\nabla u(x) \in T_{\phi(x)} \mathcal{M}$, so $P(x) \nabla u(x) = 0$ a.e. in $\Omega$, i.e. $u \in H^1_s(\Omega)$.

It remains to prove that if $v \in L^2_{\text{loc}}(\mathcal{G})$ and $\mu^{1/2} v \in L^2(\mathcal{G})$, then there exists a sequence $(u_m)_{m \in \mathbb{N}}$ in $H^1_s(\Omega)$ such that $u_m \to u$ in $L^2(\Omega)$ as $m \to \infty$, where $u = v \circ \phi$. Choose a sequence $(v_m)_{m \in \mathbb{N}}$ in $C^1_0(\mathcal{G})$, such that $v_m \to v$ in $L^2_{\text{loc}}(\mathcal{G})$. Set $u_m := v_m \circ \phi$. Then $u_m \in H^1_s(\Omega)$. Furthermore,

$$\int_\Omega J_0(x)|u(x) - u_m(x)|^2 \, \mathrm{d}x = \int_\mathcal{G} \mu(p)|v(p) - v_m(p)|^2 \, \mathrm{d}\mathcal{H}^k(p) \to 0 \quad \text{as } m \to \infty.$$

Since $\inf_{x \in \Omega} J_0(x) > 0$, the the proof is complete. $\square$

*Proof of Theorem 3.4.* This is an easy consequence of Theorem 3.3, Proposition 4.4 and the coarea formula. Trivial details are omitted. $\square$

We finally give a

*Proof of Theorem 3.7.* Set $n := k + 1$ and let $(\lambda_j)_{j \in \mathbb{N}}$ be the repeated sequence of the eigenvalues of the operator $-\Delta_{\mathbb{S}^{n-1}(r)}$. Moreover, for $\nu \in \mathbb{N}_0$, let $\overline{\lambda}_\nu$ denote the $\nu$-th distinct eigenvalue of $-\Delta_{\mathbb{S}^{n-1}(r)}$. It is well known (see e.g. [3]) that

$$\overline{\lambda}_\nu = r^{-2} \nu(\nu + n - 2), \quad \text{for } \nu \in \mathbb{N}.$$

The multiplicity of $\overline{\lambda}_\nu$ is

$$\binom{\nu + n - 1}{\nu} - \binom{\nu + n - 2}{\nu - 1}$$

and the eigenspace of $\overline{\lambda}_\nu$ is precisely the space of all homogeneous harmonic polynomials on $\mathbb{R}^n$ of degree $\nu$, restricted to $\mathbb{S}^{n-1}(r)$. So we can find arbitrarily large gaps in the spectrum of $\Delta_{\mathbb{S}^{n-1}(r)}$. In particular, we have that

$$\lim_{\nu \to \infty} \frac{\overline{\lambda}_{\nu+1} - \overline{\lambda}_\nu}{\overline{\lambda}_\nu^{1/2}} = \frac{2}{r}, \tag{20}$$

and hence

$$\limsup_{j \to \infty} \frac{\lambda_{j+1} - \lambda_j}{\lambda_j^{1/2}} = \frac{2}{r}. \tag{21}$$

Since $\mathbb{S}^{n-1}(r)$ is a manifold without boundary, then $D(A_\mu) = H^2(\mathbb{S}^{n-1}(r))$. It follows that

$$A_\mu u = -(1/\mu)\operatorname{div}(\mu \nabla u) = -\Delta_{\mathbb{S}^{n-1}(r)} u - \langle (1/\mu)\nabla \mu, \nabla u \rangle.$$

This means that $A_\mu$ is a relatively bounded perturbation of $-\Delta_{\mathbb{S}^{n-1}(r)}$. More precisely, set $A := -\Delta_{\mathbb{S}^{n-1}(r)}$ and, for $u \in H^1(\mathbb{S}^{n-1}(r))$, set

$$B_\mu u := -(1/\mu)\langle \nabla \mu, \nabla u(p) \rangle,$$

so $A_\mu = A + B_\mu$. For $u \in H^2(\mathbb{S}^{n-1}(r))$, we have that

$$|B_\mu u|_{L^2}^2 = \int_{\mathbb{S}^{n-1}(r)} |\langle \mu^{-1}\nabla \mu, \nabla u \rangle|^2 \, d\mathcal{H}^{n-1}$$

$$\leq C_\mu^2 \int_{\mathbb{S}^{n-1}(r)} \langle \nabla u, \nabla u \rangle \, d\mathcal{H}^{n-1} = C_\mu^2 \int_{\mathbb{S}^{n-1}(r)} u\, Au \, d\mathcal{H}^{n-1} \leq C_\mu^2 |u|_{L^2}|Au|_{L^2}.$$

It follows that, whenever $\delta > 0$, we have

$$|B_\mu u|_{L^2} \leq \delta |Au|_{L^2} + \frac{C_\mu^2}{4\delta}|u|_{L^2} \quad \text{for all } u \in D(A). \tag{22}$$

Now let $\lambda > 0$ and let $d(\lambda)$ be the distance of $\lambda$ from the spectrum of $A$. Assume that $\lambda I - A$ is invertible. Write $L^2 := L^2(\mathbb{S}^{n-1}(r))$. It is well known (see e.g. Theorem 3.17 in [9]) that a sufficient condition for $\lambda I - (A + B_\mu)$ being invertible is

$$|B_\mu(\lambda I - A)^{-1}|_{\mathcal{L}(L^2,L^2)} < 1.$$

In view of (22), for every $\delta > 0$ we have

$$|B_\mu(\lambda I - A)^{-1}|_{\mathcal{L}(L^2,L^2)} \leq \delta |A(\lambda I - A)^{-1}|_{\mathcal{L}(L^2,L^2)} + \frac{C_\mu^2}{4\delta}|(\lambda I - A)^{-1}|_{\mathcal{L}(L^2,L^2)}.$$

Observe that, since $A$ is self-adjoint,

$$|(\lambda I - A)^{-1}|_{\mathcal{L}(L^2,L^2)} = \sup_{\nu \in \mathbb{N}} |\lambda - \overline{\lambda}_\nu|^{-1} \leq d(\lambda)^{-1}$$

and

$$|A(\lambda I - A)^{-1}|_{\mathcal{L}(L^2,L^2)} = \sup_{\nu \in \mathbb{N}} |\overline{\lambda}_\nu||\lambda - \overline{\lambda}_\nu|^{-1}$$

$$\leq \sup_{\nu \in \mathbb{N}}(1 + \lambda|\lambda - \overline{\lambda}_\nu|^{-1}) \leq 1 + \lambda d(\lambda)^{-1}.$$

It follows that

$$|B_\mu(\lambda I - A)^{-1}|_{\mathcal{L}(L^2,L^2)} \leq \delta(1 + \lambda d(\lambda)^{-1}) + \frac{C_\mu^2}{4\delta}d(\lambda)^{-1}.$$

So a sufficient condition for $\lambda I - (A + B_\mu)$ being invertible is

$$\delta(d(\lambda) + \lambda) + \frac{C_\mu^2}{4\delta} < d(\lambda)$$

or equivalently

$$(23) \qquad \delta\lambda + \frac{C_\mu^2}{4\delta} < (1-\delta)d(\lambda) \quad \text{for some } \delta, \ 0 < \delta < 1$$

Using our assumption on $C_\mu$ we see that (23) is satisfied (and so $\lambda I - A_\mu$ is invertible) whenever

$$(24) \qquad \lambda > 1/(4r)^2 \quad \text{and} \quad d(\lambda) > \frac{1}{2r}\lambda^{1/2}.$$

(To see this, just set $\delta := (8r)^{-1}\lambda^{1/2}$.) Now let $\nu > 1$ be fixed. Then $\nu(\nu + n - 2) > 1/4$, so $\overline{\lambda}_\nu > 1/(4r)^2$. Let us consider the interval $]\overline{\lambda}_\nu, \overline{\lambda}_{\nu+1}[$. If $\lambda \in ]\overline{\lambda}_\nu, \overline{\lambda}_{\nu+1}[$, then, in view of (24), $\lambda I - A_\mu$ is invertible provided

$$\lambda - \overline{\lambda}_\nu > \frac{1}{2r}\lambda^{1/2} \quad \text{and} \quad \overline{\lambda}_{\nu+1} - \lambda > \frac{1}{2r}\lambda^{1/2}.$$

Set $\xi := \lambda - \overline{\lambda}_\nu$ and $\eta := \overline{\lambda}_{\nu+1} - \lambda$. Thus $\lambda I - A_\mu$ is invertible provided

$$4r^2\xi^2 - \xi - \overline{\lambda}_\nu > 0 \quad \text{and} \quad 4r^2\eta^2 + \eta - \overline{\lambda}_{\nu+1} > 0.$$

By solving these inequalities for $\xi$ and $\eta > 0$, we obtain the conditions

$$\xi > \xi_\nu := \frac{1}{8r^2} + \left(\frac{1}{64r^4} + \frac{\overline{\lambda}_\nu}{4r^2}\right)^{1/2}$$

and

$$\eta > \eta_{\nu+1} := -\frac{1}{8r^2} + \left(\frac{1}{64r^4} + \frac{\overline{\lambda}_{\nu+1}}{4r^2}\right)^{1/2}.$$

It follows that, if $\overline{\lambda}_\nu + \xi_\nu < \overline{\lambda}_{\nu+1} - \eta_{\nu+1}$, then the interval

$$I_\nu := ]\overline{\lambda}_\nu + \xi_\nu, \overline{\lambda}_{\nu+1} - \eta_{\nu+1}[$$

is contained in the resolvent set of $A_\mu$. So let us compute

$$(\overline{\lambda}_{\nu+1} - \eta_{\nu+1}) - (\overline{\lambda}_\nu + \xi_\nu)$$

$$= \overline{\lambda}_{\nu+1} - \overline{\lambda}_\nu - \left(\left(\frac{1}{64r^4} + \frac{\overline{\lambda}_{\nu+1}}{4r^2}\right)^{1/2} + \left(\frac{1}{64r^4} + \frac{\overline{\lambda}_\nu}{4r^2}\right)^{1/2}\right)$$

$$= \overline{\lambda}_{\nu+1} - \overline{\lambda}_\nu - \frac{1}{4r^2}(\overline{\lambda}_{\nu+1} - \overline{\lambda}_\nu)\left(\left(\frac{1}{64r^4} + \frac{\overline{\lambda}_{\nu+1}}{4r^2}\right)^{1/2} - \left(\frac{1}{64r^4} + \frac{\overline{\lambda}_\nu}{4r^2}\right)^{1/2}\right)^{-1}.$$

Substituting the explicit expression $\overline{\lambda}_\nu = r^{-2}\nu(\nu + n - 2)$, a straightforward computation shows that

$$\lim_{\nu \to \infty}\left(\left(\frac{1}{64r^4} + \frac{\overline{\lambda}_{\nu+1}}{4r^2}\right)^{1/2} - \left(\frac{1}{64r^4} + \frac{\overline{\lambda}_\nu}{4r^2}\right)^{1/2}\right) = \frac{1}{2r^2}.$$

It follows that there is a $\nu_0 \in \mathbb{N}$ such that for all $\nu \geq \nu_0$,

$$(25) \qquad (\overline{\lambda}_{\nu+1} - \eta_{\nu+1}) - (\overline{\lambda}_\nu + \xi_\nu) \geq \frac{1}{3}(\overline{\lambda}_{\nu+1} - \overline{\lambda}_\nu).$$

In particular, for such $\nu$, the interval $I_\nu$ is nonempty.

Let
$$\lambda_1^\mu \leq \lambda_2^\mu \leq \lambda_3^\mu \leq \ldots$$
be the repeated sequence of the eigenvalues of $A_\mu$.

For $\nu \geq \nu_0 + 1$, set
$$J_\nu := ]\overline{\lambda}_\nu - \eta_\nu, \overline{\lambda}_\nu + \xi_\nu[.$$

Define $Z$ to be the set of all $\nu \geq \nu_0 + 1$ such that $J_\nu$ has nonempty intersection with the spectrum of $A_\mu$. It follows that $Z$ has infinitely many elements. For $\nu \in Z$ set
$$j_\nu := \max\{\, j \in \mathbb{N} \mid \lambda_j^\mu \in J_\nu \,\}.$$

$j_\nu$ is well defined since $\lambda_j^\mu \to \infty$ as $j \to \infty$. Now $\lambda_{j_\nu+1}^\mu \geq \lambda_{j_\nu}^\mu$ and so, by the definition of $j_\nu$ and the fact that there are no eigenvalues of $A_\mu$ lying in $I_\nu$, it follows that $\lambda_{j_\nu+1}^\mu \geq \overline{\lambda}_{\nu+1} - \eta_{\nu+1}$. Therefore we have that
$$\frac{\lambda_{j_\nu+1}^\mu - \lambda_{j_\nu}^\mu}{(\lambda_{j_\nu}^\mu)^{1/2}} \geq \frac{(\overline{\lambda}_{\nu+1} - \eta_{\nu+1}) - (\overline{\lambda}_\nu + \xi_\nu)}{(\overline{\lambda}_\nu + \xi_\nu)^{1/2}}.$$

Thus, by (25), we see that
$$\frac{\lambda_{j_\nu+1}^\mu - \lambda_{j_\nu}^\mu}{(\lambda_{j_\nu}^\mu)^{1/2}} \geq \frac{1}{3} \frac{\overline{\lambda}_{\nu+1} - \overline{\lambda}_\nu}{\overline{\lambda}_\nu^{1/2}} \frac{\overline{\lambda}_\nu^{1/2}}{(\overline{\lambda}_\nu + \xi_\nu)^{1/2}}.$$

Since
$$\lim_{\nu \to \infty} \frac{\overline{\lambda}_\nu^{1/2}}{(\overline{\lambda}_\nu + \xi_\nu)^{1/2}} = 1,$$
we obtain, in view of (20), that
$$\limsup_{\nu \to \infty} \frac{\lambda_{j_\nu+1}^\mu - \lambda_{j_\nu}^\mu}{(\lambda_{j_\nu}^\mu)^{1/2}} \geq \frac{2}{3r}$$
and therefore
$$(26) \qquad \limsup_{j \to \infty} \frac{\lambda_{j+1}^\mu - \lambda_j^\mu}{(\lambda_j^\mu)^{1/2}} \geq \frac{2}{3r} > 0.$$

Now, in view of Theorem 3.4, the repeated eigenvalue sequences of the limit operator $A_0$ and the operator $A_\mu$ are the same. The proof is complete. $\square$

Università degli Studi di Trieste, Dipartimento di Scienze Matematiche, Via Valerio, 12/b, 34100 Trieste, Italy

*E-mail address*: `prizzi@mathsun1.univ.trieste.it`

Universität Rostock, Fachbereich Mathematik, Universitätsplatz 1, 18055 Rostock, Germany

*E-mail address*: `krzysztof.rybakowski@mathematik.uni-rostock.de`